\documentclass[letterpaper,11pt,onecolumn]{article}

\usepackage[margin=1in]{geometry}

\usepackage[T1]{fontenc}
\usepackage[latin1]{inputenc}
\usepackage[english]{babel}

\usepackage{amsmath}
\usepackage{amsbsy}
\usepackage{amsfonts}
\usepackage{amssymb}
\usepackage{bm}
\usepackage{booktabs}
\usepackage{csquotes}
\usepackage{diagbox}
\usepackage{enumerate}
\usepackage{extarrows}
\usepackage{float}
\usepackage{graphicx}
\usepackage{mathrsfs}
\usepackage{pbox}
\usepackage{pgf}
\usepackage{xcolor}
\usepackage{xfrac}

\usepackage[font={footnotesize}]{caption}
\usepackage{subcaption}

\usepackage{placeins} % FloatBarrier

\usepackage[pdftex]{hyperref}
\hypersetup{colorlinks=false}

\usepackage{array}
\newcolumntype{L}[1]{>{\raggedright\let\newline\\\arraybackslash\hspace{0pt}}p{#1}}
\newcolumntype{C}[1]{>{\centering\let\newline\\\arraybackslash\hspace{0pt}}p{#1}}
\newcolumntype{R}[1]{>{\raggedleft\let\newline\\\arraybackslash\hspace{0pt}}p{#1}}

\newcommand{\idmatrix}{\textup{\uppercase\expandafter{\romannumeral 1}}}

\makeatletter
\newcommand{\dashto}[1][2pt]{
  \settowidth{\@tempdima}{${}\rightarrow{}$}
  \makebox[\@tempdima]{${}\rightarrow{}$}% typeset arrow
  \makebox[-\@tempdima]{\hspace{-0.1\@tempdima}\color{white}\rule[0.5ex]{#1}{1pt}}% typeset overlay
  \makebox[\@tempdima]{}% advance appropriate horizontal distance
  }
\makeatother

\def\gobench{\texttt{GO\_3D\_OBS}}

\newcommand{\fine}[1]{\texttt{fine(#1)}}
\newcommand{\coarse}[1]{\texttt{coarse(#1)}}
  \definecolor{ICES}{RGB}{94, 156, 174}
  \definecolor{ORANGE}{RGB}{191, 87, 0}
  \definecolor{RED}{RGB}{190, 30, 49}
  \definecolor{SUN}{RGB}{227, 81, 51}
  \definecolor{GREEN}{RGB}{0, 171, 86}
  \definecolor{BLUE}{RGB}{11, 78, 179}
  \definecolor{BROWN}{RGB}{122, 80, 40}
  \definecolor{GREY}{RGB}{50, 50, 50}
  \definecolor{TEAL}{RGB}{0, 160, 176}

\newcommand{\bb}[1]{\mathbb{#1}}
\newcommand{\mc}[1]{\mathcal{#1}}

\newcommand{\Gammah}{\Gamma_{\hskip -1pt h}}

\def\gradh{\nabla_{\hspace{-1pt}h}}

\def\tcurl{\text{curl}}
\def\tdiv{\text{div}}

\def\hH1{H^1(\mathcal{T}_h)}

\def\H1{H^1(\Omega)}

\numberwithin{equation}{section}

\graphicspath{{figures/}}

\usepackage[backend=bibtex,
	giveninits=true,
	maxbibnames=4,
	date=year,
	doi=false,
	url=false,
	isbn=false,
	eprint=false]{biblatex}
\begin{document}

\baselineskip=16pt
\parskip=4pt

%=====================================================%
% Title and authors
\begin{center}
{\Large {\bf  
Scalable DPG Multigrid Solver for Helmholtz Problems:\\
A Study on Convergence}}\\[8pt]
{\Large
	Jacob Badger$^a$,
	Stefan Henneking$^a$,
	Socratis Petrides$^b$,
	Leszek Demkowicz$^a$
	\\[12pt]}
{\large
$^a$Oden Institute, The University of Texas at Austin \\
$^b$Lawrence Livermore National Laboratory
}
\end{center}

\vspace{10pt}

\begin{center}
{\large \textbf{Dedicated to Leszek Demkowicz on the occasion of his 70$^{\textbf{th}}$ birthday}}
\end{center}

\vspace{10pt}

%=====================================================%
% Abstract
\begin{abstract}
%
%!TEX root = ../paper.tex
%

This paper presents a scalable multigrid preconditioner targeting large-scale systems arising from discontinuous Petrov--Galerkin (DPG) discretizations of high-frequency wave operators. This work is built on previously developed multigrid preconditioning techniques  of Petrides and Demkowicz (Comput.~Math.~Appl.~87 (2021) pp.~12--26) and extends the convergence results from $\mathcal{O}(10^7)$ degrees of freedom (DOFs) to $\mathcal{O}(10^9)$ DOFs using a new scalable parallel MPI/OpenMP implementation. Novel contributions of this paper include an alternative definition of coarse-grid systems based on restriction of fine-grid operators, yielding superior convergence results. In the uniform refinement setting, a detailed convergence study is provided, demonstrating $h$ and $p$ robust convergence and linear dependence with respect to the wave frequency. The paper concludes with numerical results on $hp$-adaptive simulations including a large-scale seismic modeling benchmark problem with high material contrast.

\vspace{20pt}


\vspace{20pt}

\end{abstract}

%=====================================================%
% Acknowledgments
\subsection*{Acknowledgments}

We thank Matteo Croci for helpful discussions.
J. Badger, S. Henneking, and L. Demkowicz were supported with AFOSR grant FA9550-19-1-0237 and NSF award 2103524.
This work was performed under the auspices of the U.S. Department of Energy by Lawrence Livermore National Laboratory under Contract DE-AC52-07NA27344, LLNL-JRNL-847116.

%=====================================================%
% Sections
%
%!TEX root = ../paper.tex
%

\section{Introduction}
\label{sec:intro}

%=====================================================%
\subsection{Background}

Wave propagation problems arise in a number of contexts including natural resource exploration, medical imaging, and nuclear fusion research, to name a few. However, developing accurate and efficient numerical algorithms for the solution of time-harmonic wave propagation problems is a notoriously difficult problem. While traditional finite element methods (FEM) can deliver high-accuracy and optimal discretizations, their efficacy for wave operators deteriorates for two main reasons. First, they suffer from stability issues unless very fine meshes are used to resolve the propagating wave. In the high-frequency regime this results in prohibitively expensive problems. The lack of preasymptotic discrete stability also makes mesh adaptivity techniques unreliable and inefficient. Second, the linear system is highly indefinite and, consequently, standard iterative solution schemes break down \cite{ernst2012difficult}. Current leading-edge preconditioning techniques for wave operators, such as multigrid methods \cite{hiptmair1998multigrid, stolk2014multigrid, gopalakrishnan2004analysis}, domain decomposition methods with special transmission conditions \cite{cai1992domain, gopalakrishnan2003overlapping, stolk2013rapidly, kim2015optimized}, stabilized methods based on artificial absorption \cite{absorption_1, absorption_2}, shifted Laplacian \cite{sheikh2013convergence} and sweeping preconditioners \cite{engquist2011sweeping1, engquist2011sweeping2, chen2013source, leng2020diagonal, taus2020sweeps} are very promising but they lose their efficiency in heterogeneous media and can be difficult to apply in complex geometries \cite{gander2019class, ErlanggaYogiA}.

An alternative approach instead employs \emph{minimum-residual} discretization methodologies which, by construction, produce positive-definite discrete systems and may therefore be amenable to more standard preconditioning techniques \cite{Gopalakrishnan2015,petrides2017adaptive,petrides2021adaptive,petrides2019phd}. Indeed, popularization of the first-order system least-squares methodology (FOSLS) \cite{cai1997first, lee2000first}, and other least-squares methodologies \cite{de2004least}, was driven by the applicability of geometric and algebraic multigrid methods to otherwise indefinite problems. However, for wave propagation problems, FOSLS is known to be highly dissipative \cite{gopalakrishnan2014dispersive} and thus not competitive in the high-frequency regime. This work discusses a multigrid solver based on a minimum-residual discretization obtained by the discontinuous Petrov--Galerkin (DPG) method with Optimal Test Functions \cite{demkowicz2017dpg} applied to the ultraweak variational formulation.

%=====================================================%
\subsection{DPG-MG Solver}

The DPG FE methodology of Demkowicz and Gopalakrishnan \cite{demkowicz2010part1, demkowicz2011part2, demkowicz2012part3} is a non-standard least-squares method with several attractive properties: mesh-independent stability, a built-in error indicator, and applicability to a number of variational formulations with different functional settings. A special case of the DPG method is the well-established FOSLS method in which the residual is minimized in the $L^2$ test norm. As mentioned previously, however, other formulations are preferable in the context of wave propagation. Among the various DPG formulations, the so-called \emph{ultraweak} variational formulation has proved to be superior: it is less dissipative than other DPG formulations \cite{petrides2019phd}, with dispersion error roughly commensurate to Galerkin discretizations \cite{gopalakrishnan2014dispersive}, and it has been shown to solve problems with many wavelengths accurately by countering the pollution error through a modest increase in the order of discretization \cite{henneking2021pollution}. These properties were leveraged by Petrides and Demkowicz in \cite{petrides2021adaptive} to define an $hp$-adaptive multilevel preconditioner for DPG wave propagation problems discretized with conforming elements of the exact-sequence energy spaces \cite{demkowicz2020fem}.

Similar to hybridizable methods, the DPG methodology introduces additional trace degrees of freedom (DOFs) on the mesh skeleton resulting from testing with larger discontinuous\footnote{The `D' in the DPG name} (``broken'') test spaces \cite{carstensen2016breaking}. In the case of high-order discretizations, statically condensing all interior DOFs onto the mesh skeleton results in a smaller global system and enables more coherent implementation for field and trace variables. The DPG multigrid solver (DPG-MG) is defined on this condensed global system of trace degrees of freedom. Constructing suitable prolongation operators for the condensed system is complicated by the fact that fine-grid DOFs resulting from $h$-refinement have no natural coarse-grid representatives; this is a challenge shared by hybridizable methods \cite{rhebergen2018preconditioning, rhebergen2022preconditioning}. The construction of a stable prolongation operator between such non-nested condensed systems for general DPG problems\footnote{discretized with exact-sequence energy spaces} was one of the major contributions in the original DPG-MG work by Petrides and Demkowicz \cite{petrides2017adaptive, petrides2019phd} and it will be outlined later in Section~\ref{sec:DPG-MG}. 

%=====================================================%
\subsection{Direction and Outline}

Based on the initial implementation by Petrides and Demkowicz \cite{petrides2019phd}, we have developed a scalable hybrid MPI/OpenMP implementation of the DPG-MG solver. The parallel implementation and its scaling characteristics will be detailed in a forthcoming publication; the present work instead leverages our performant implementation to study the convergence properties of the DPG-MG solver under uniform $h$, uniform $p$, and $hp$-adaptive refinements. This work is intended to elucidate scaling characteristics of the DPG-MG solver and identify aspects of the current construction which may be improved. 

The remainder of the paper is structured as follows: Section~\ref{sec:problem} defines the ultraweak acoustics model and DPG discretization used throughout this work. Section~\ref{sec:DPG-MG} reviews the construction of the DPG-MG solver. In Section~\ref{sec:studies}, a number of convergence studies are performed for a model problem with manufactured solution. The applicability of the solver to state-of-the-art computational challenges is demonstrated in Section~\ref{sec:application}, using the \gobench{} model \cite{go_3d_obs}; a challenging benchmark problem for evaluating next-generation algorithms in seismic modeling.
We conclude in Section~\ref{sec:conclusion} with a discussion of findings and future work.

%
%!TEX root = ../paper.tex
%

\section{Ultraweak DPG for Helmholtz}
\label{sec:problem} 

The DPG method constructs automatically stable discretizations of well-posed variational formulations, inheriting stability from the continuous problem. It achieves this by computing optimal test functions \cite{demkowicz2011part2} that realize the supremum in the discrete inf--sup condition \cite{babuska1971error}. The unique space of these specially-selected test functions is called the optimal test space. In practice, this space is approximated by inverting the global Riesz map over an enriched, discontinuous test space \cite{gopala2014practical}. The numerical computations in this paper employ a uniform increase of the polynomial order by $1$ for the enrichment of the test space. The discontinuous (broken) nature of the test space enables the element-local computation of optimal test functions. However, this breaking of the test space results in additional (trace) unknowns defined on the mesh skeleton \cite{carstensen2016breaking}.

\paragraph{Notation and energy spaces.}

We briefly introduce some notation and the energy spaces used throughout this work. Consider a bounded domain $\Omega \subset \bb{R}^3$ with Lipschitz boundary $\Gamma \equiv \partial \Omega$. The $L^2$-inner product over $\Omega$ is denoted by $(\cdot,\cdot)$ and the $L^2$-norm by $\| \cdot \|$. We define the standard energy spaces
\begin{equation}
\begin{alignedat}{2}
	L^2(\Omega) &= \{ y: \Omega \rightarrow \bb{C} : \| y \| < \infty \} , \\
	H^1(\Omega) &= \{ w : \Omega \rightarrow \bb{C} : w \in L^2(\Omega), \nabla w \in (L^2(\Omega))^3 \} , \\
	H(\tdiv, \Omega) &= \{ v: \Omega \rightarrow \bb{C}^3 : v \in (L^2(\Omega))^3, \nabla \cdot v \in L^2(\Omega) \} .
\end{alignedat}
\end{equation}

In the DPG method, we use corresponding broken energy spaces for test functions which are defined as product-spaces over elements $\{ K \}_{K \in \Omega_h}$ of the finite element mesh $\Omega_h$:
\begin{equation}
\begin{alignedat}{2}
	H^1(\Omega_h) &:= \{ w : \Omega \rightarrow \bb{C} : w|_K \in H^1(K)\ \forall K \in \Omega_h \} , \\
	H(\tdiv, \Omega_h) &:= \{ v: \Omega \rightarrow \bb{C}^3 : v|_K \in H(\tdiv, K)\ \forall K \in \Omega_h \} .
\end{alignedat}
\end{equation}

Lastly, the breaking of test functions \cite{carstensen2016breaking} leads to introducing trace unknowns on the mesh skeleton $\Gammah := \{ \partial K \}_{K \in \Omega_h}$. The trace spaces are understood as element-wise traces of globally conforming functions:
\begin{equation}
\begin{alignedat}{2}
	H^{1/2}(\Gammah) &:=
	\{ \prod_{K \in \Omega_h} \gamma^K (w |_K) : w \in H^1(\Omega) \} , \\
	H^{-1/2}(\Gammah) &:=
	 \{ \prod_{K \in \Omega_h} \gamma_n^K (v |_K) : v \in H(\tdiv,\Omega) \} ,
\end{alignedat}
\end{equation}
where $\gamma^K$ and $\gamma_n^K$ are element-wise continuous and normal trace operators \cite{demkowicz2017dpg, demkowicz2018spaces}.

\paragraph{Helmholtz problem.}
This work considers the first-order mixed form of time-harmonic linear acoustics with inhomogeneous impedance boundary condition (BC). In operator form, the equations are given by
\begin{equation}
\begin{alignedat}{3}
	i \omega p + \nabla \cdot u &= 0 \quad &&\text{in } \Omega, \\
	i \omega u + \nabla p &= 0 \quad &&\text{in } \Omega, \\
	Z^{-1} p - u_n &= u_0 \quad &&\text{on } \Gamma ,
\end{alignedat}
\label{eq:helmholtz}
\end{equation}
where $p$ is pressure, $u$ is velocity, $\omega$ is the angular wave frequency, and $i = \sqrt{-1}$; in the impedance BC, $Z$ is acoustic wave impedance, $u_n := u \cdot n$ is the flux in outward normal direction $n$, and $u_0$ is an impedance load.

\paragraph{Broken ultraweak formulation.}

Let $(\mathcal{U}, \hat{\mathcal{U}})$ and $\mathcal{V}$ be the trial and test space, respectively, and let $\mathcal{V}'$ be the space of antilinear functionals on $\mathcal{V}$. The DPG formulation of the Helmholtz problem is defined by a variational formulation of the form: Given $l \in \mathcal{V}'$, find $\mathfrak u \in \mathcal{U}$ and $\hat{\mathfrak u} \in \hat{\mathcal{U}}$ that satisfy
\begin{equation}
	b(\mathfrak u, \mathfrak v) + \hat{b}(\hat{\mathfrak u},\mathfrak v) = l(\mathfrak v), \quad \mathfrak v \in \mathcal{V} ,
	\label{eq:variational-formulation}
\end{equation}
where $b$ and $\hat b$ are sesquilinear forms on $\mathcal{U} \times \mathcal{V}$ and $\hat{\mathcal{U}} \times \mathcal{V}$, respectively.

We refer to \cite{petrides2019phd, demkowicz2020fem} for a thorough derivation of the ultraweak Helmholtz formulation. The broken ultraweak formulation, given by (\ref{eq:variational-formulation}), is defined by the following group variables and forms:
\begin{equation}
\begin{split}
	\mathfrak u = &\ (u, p) \in (L^2(\Omega))^3 \times L^2(\Omega) , \\
	\hat{\mathfrak u} =&\ (\hat u_n, \hat{p}) \in H^{-\frac{1}{2}}(\Gammah) \times H^{\frac{1}{2}}(\Gammah) : 
	Z^{-1} \hat p - \hat u_n = 0 \text{ on } \Gamma , \\
	\mathfrak v =&\ (q, v) \in H^1(\Omega_h) \times H(\tdiv, \Omega_h) :
	Z^{-1} q + v_n = 0 \text{ on } \Gamma , \\[5pt]
	%%%%%%
	b(\mathfrak u, \mathfrak v) =&\ (i \omega p, q) - (u, \gradh q) + (i \omega u, v) - (p, \gradh \cdot v) \\
	\hat{b}(\hat{\mathfrak u}, \mathfrak v) =&\ 
	\langle \hat u_n, q \rangle_{\Gammah} +
	\langle \hat p, v_n \rangle_{\Gammah} , \\
	l(\mathfrak v) =&\ \langle u_0, q \rangle_{\Gamma} .
\end{split}
\end{equation}
Note that the test functions are assumed to satisfy $Z^{-1} q + v_n = 0 \text{ on } \Gamma$ in order to build in the impedance BC, and that the impedance BC implicitly implies additional regularity of the velocities on boundary $\Gamma$. For the load to be well-defined, we can assume $u_0 \in H^{-1/2}(\Gamma)$ in which case $\langle u_0, q \rangle_{\Gamma}$ can be understood in the sense of duality pairing; another option is to assume $u_0 \in L^2(\Gamma)$ interpreting $\langle u_0, q \rangle_{\Gamma}$ in the $L^2$-sense; see \cite{demkowicz2020fem} for further discussion on the regularity issue for impedance BCs.

The additional unknowns $\hat u_n$, $\hat p$ describe the normal velocity (i.e., flux) and the pressure on element boundaries on the mesh skeleton $\Gammah$; $\hat u_n$ is discretized as the normal trace of $H(\tdiv)$-conforming elements, and $\hat p$ as the continuous trace of $H^1$-conforming elements. The broken test space is equipped with the adjoint graph norm \cite{petrides2021adaptive, demkowicz2017dpg}:
\begin{equation}
	\| \mathfrak v \|_{\mathcal V}^2 :=
	\| A_h^* \mathfrak v \|^2 + \alpha \| \mathfrak v \|^2 ,
	\label{eq:Gnorm}
\end{equation}
where $A_h^* \mathfrak v = -(i \omega q + \gradh \cdot v, \gradh q + i \omega v)$, and $\alpha$ is a scaling constant. Throughout this work, numerical results are computed with $\alpha=1$.

%
%!TEX root = ../paper.tex
%

\section{DPG-MG Solver}
\label{sec:DPG-MG}

DPG-MG is a multigrid-preconditioned conjugate gradient solver, with the multilevel preconditioner defined on a hierarchy of meshes produced through refinement as defined in \cite{petrides2021adaptive}. Mesh-independent stability of the DPG methodology implies that the DPG-MG solver can be initialized on arbitrarily coarse initial meshes. Once the solution is obtained to a sufficient accuracy on the current mesh, the DPG error indicator is used to produce a set of refinements to define the next mesh. Because the solution on intermediate meshes is needed only to sufficient accuracy to produce the next mesh, optimal $hp$-adaptive meshes can be produced with relatively few iterations and at little cost.

\paragraph{Prolongation.} As indicated in Section \ref{sec:intro}, the DPG-MG solver is defined on trace DOFs located on the mesh skeleton; $h$-refinements produce fine-grid edges and faces that do not coincide with the previous-grid skeleton and thus have no natural representatives on the previous mesh. To ameliorate this, Petrides and Demkowicz introduced a two-stage prolongation \cite{petrides2017adaptive}. In the first stage of the restriction, fine-grid DOFs not supported on the previous mesh skeleton are statically condensed; the resulting mesh is called the \emph{macro} grid. In the second stage of the restriction, a natural inclusion operator---that expresses coarse-grid basis functions as a linear combination of macro-grid basis functions---is applied on edges and faces of the previous mesh. The natural inclusion operator is constructed via recursive application of \emph{constrained approximation} techniques \cite{hpbook2}.

\paragraph{V-cycle and smoother.}

\begin{figure}[htb]
	\centering
	\includegraphics[width=\textwidth]{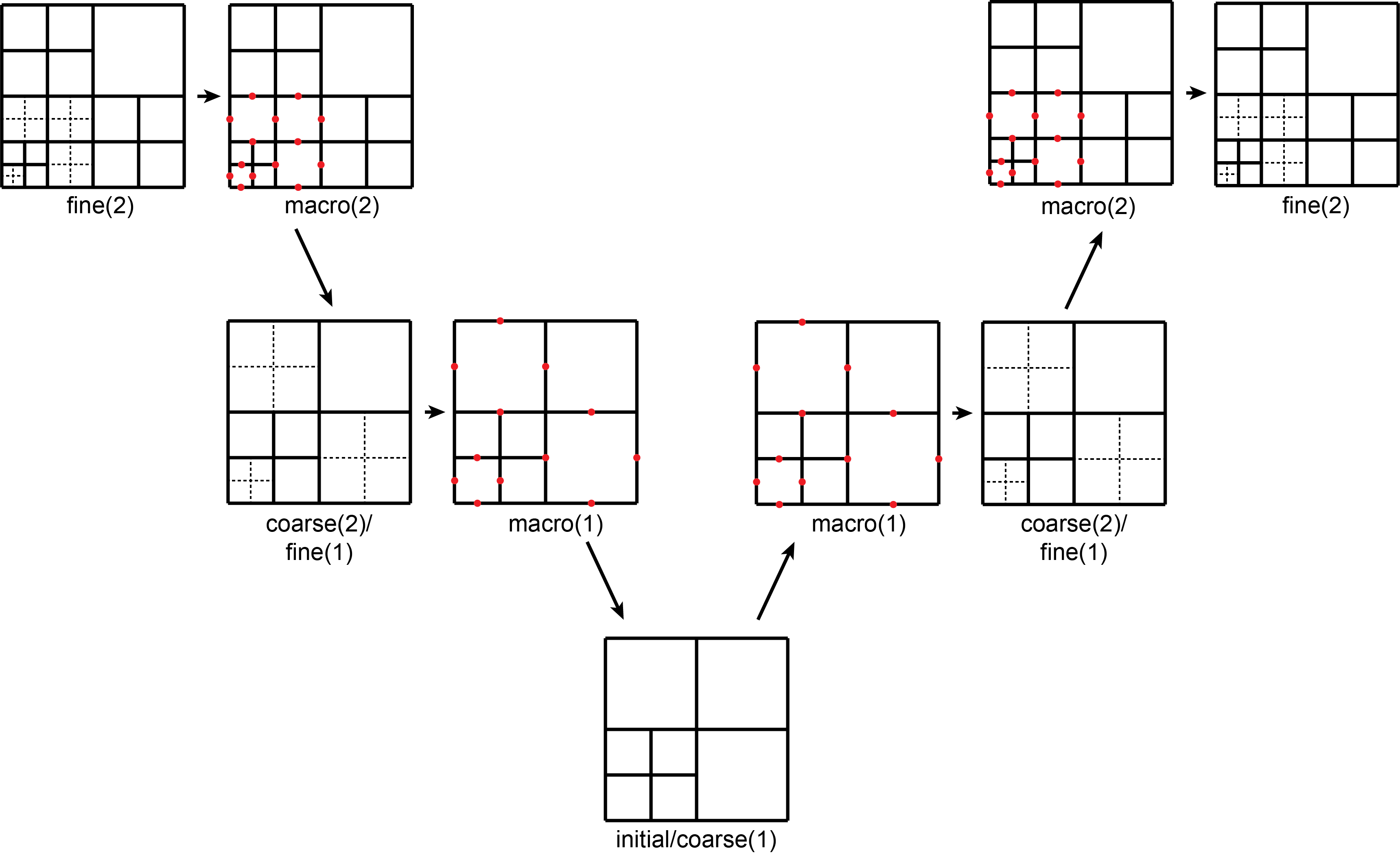}
	\caption{DPG-MG V-cycle. The prolongation operator is defined by first condensing all fine-grid degrees of freedom onto the coarse-grid skeleton; then applying the natural inclusion operator on macro-grid edges and, in 3D, faces.}
	\label{fig:v-cycle}
\end{figure}

We briefly illustrate the DPG-MG V-cycle and introduce terminology used to refer to the various meshes involved in the algorithm. First, we depart somewhat from convention and refrain from designating the \emph{initial} mesh the coarse grid (we refer to it simply as the initial mesh). Therefore, terms like \emph{coarse}-grid vertex patches and \emph{coarse}-grid stiffness matrices are not necessarily related to the initial mesh. Instead, we reserve the terms \emph{coarse} and \emph{fine} grid to indicate the coarser and finer meshes, respectively, in a pair of consecutive meshes. In particular, \fine{i} refers to the grid produced by the $i$-th refinement of the initial mesh with \coarse{i} being its coarse counterpart. This naming convention is illustrated along with the V-cycle in Fig.~\ref{fig:v-cycle}. Note that in the MPI-distributed implementation, discussed below, \coarse{i} and \fine{i-1} correspond to two different partitions of the same mesh.

The DPG-MG solver applies both conjugate gradient iterations and smoothing on the macro grid. Smoothing is performed using an additive Schwarz (overlapping block-Jacobi) smoother with blocks defined as macro-grid DOFs supported on coarse-grid vertex patches. Smoothing patches on \fine{i} are constructed using \coarse{i} vertex patches. Use of vertex patches can lead to large smoothing blocks, especially in the case of high-order discretizations, but avoids additional complexities for vector-valued variables set in $H(\tcurl)$ and $H(\tdiv)$ \cite{arnold2000multigrid}. Thus, the DPG-MG solver is applicable, without modification, to any well-posed DPG problem discretized with exact-sequence energy spaces. Alternative definitions of smoothing patches lead to smaller patch sizes for vector-valued variables \cite{hiptmair1997multigrid, hiptmair1998multigrid}.

\paragraph{Alternative construction of coarse-grid operators.}

Coarse-grid operators can be constructed either by direct matrix assembly, or by restricting fine-grid operators. Note that because of the use of trace spaces and the \emph{on-the-fly} computation of the DPG optimal test space, these two approaches are not equivalent (indeed, in a non-nested multigrid iteration, as is the case here, these approaches are generally not equivalent). The DPG-MG solver assembles and solves the system on the current mesh, then refines to define the next mesh. Thus, for linear problems, direct assembly of coarse-grid operators can be accomplished by simply storing the current-grid system before refinement; this was the approach taken in the original DPG-MG implementation \cite{petrides2021adaptive}. However, we observed that for high-frequency problems, convergence rapidly deteriorated with increasing frequency. As will be demonstrated in Section~\ref{sec:studies}, computing coarse-grid stiffness matrices as restrictions of fine-grid stiffness matrices restores the expected convergence.  We are working to develop a rigorous understanding of how these two coarse-grid approaches differ, and thus defer an analysis of this phenomenon to a later publication.
For now it will suffice to say that, for large frequencies (and fixed $\alpha$), the spectrum of element Gram matrices with the adjoint graph test norm (\ref{eq:Gnorm}) changes dramatically as a mesh transitions from preasymptotic to asymptotic regime; this in turn may impart vastly different scales to fine- and coarse-grid systems and, when not accounted for in prolongation, may cause the coarse-grid correction to become unstable. 
Finally, we note that computing coarse-grid operators as restrictions of fine-grid operators is relatively inexpensive compared to assembly of the fine-grid DPG system, typically requiring between $1\%$ and $10\%$ of the cost of fine-grid assembly, even when sum-factorization \cite{mora2019fast,badger2020fast} is used to accelerate element assembly.

\paragraph{Scalable MPI/OpenMP implementation in $hp$3D.}

The DPG-MG solver is implemented in $hp$3D, a scalable finite element software for analysis and discretization of complex three-dimensional multiphysics applications \cite{hpUserManual, hpbook3}. $hp$3D supports a number of advanced FE technologies including exact-sequence conforming discretizations, fully anisotropic $hp$-adaptivity, and hybrid meshes with elements of ``all shapes'' (tetrahedra, hexahedra, prisms, pyramids). The code leverages hybrid MPI/OpenMP parallelism and interfaces with various scientific libraries including PETSc \cite{petsc-user-ref}, MUMPS \cite{amestoy2001mumps}, and Zoltan \cite{ZoltanOverviewArticle2002}. $hp$3D is available as an open-source code under a BSD-3 license.\footnote{\url{https://github.com/Oden-EAG/hp3d}}

The original DPG-MG implementation \cite{petrides2021adaptive} employed shared-memory parallelism via OpenMP threading for single-node computation. Memory limitations of typical compute node configurations limited the scalability of the adaptive solver to $\mathcal{O}(10^7)$ DOFs. We have extended the original solver implementation to support scalable distributed-memory computation with MPI. The approach is based on distributing solution and geometry DOFs on subdomains \cite{hpbook3,henneking2021phd}; however, the DPG-MG solver employs unique data structures and algorithms that extend $hp$3D's data structures to allow for asymmetric inclusion of ghost elements and enable efficient asynchronous communication with neighboring subdomains. To maintain satisfactory parallel efficiency on $hp$-adaptive meshes, dynamic load balancing is performed at each new grid level during the refinement process. For a fixed frequency, the distributed DPG-MG solver implementation has been shown to scale with near-linear parallel efficiency to $\mathcal{O}(10^9)$ DOFs on $\mathcal{O}(100)$ compute nodes. The details of the solver's parallel data structures, algorithms, and its scaling characteristics are not the focus of this paper and will instead be discussed in a future publication.

%
%!TEX root = ../paper.tex
%

%=====================================================%
\section{Convergence Studies}
\label{sec:studies}

In this section, we perform a number of convergence studies to investigate the robustness of the DPG-MG solver with respect to element size $h$, polynomial order $p$, and angular frequency $\omega$. Previous expositions of the DPG-MG solver considered a variety of physical problems, smoothing steps, and tolerances; illustrating the versatility of the DPG-MG solver but somewhat confounding the scaling behavior. Instead, to elucidate the convergence characteristics of the DPG-MG solver, we fix the following parameters:
\begin{itemize}
	\item Conjugate gradient iterations are terminated when the relative $\ell^2$-norm of the discrete residual has been reduced by a factor of $10^{7}$.
	\item After each refinement, the (initial) solution is reset to zero; in other words, solutions from previous grids are \emph{not} used to generate initial guesses for following grids. 
	\item A single pre- and post-smoothing step is performed on each grid level (V(1,1)-cycle), except in one case in Section \ref{sec:recomp} in which both \emph{one} and \emph{five} smoothing steps are employed to aid in comparison; this case will be noted.
	\item The initial mesh is a single element of order $p=2$; however, iterations are not reported for the initial mesh which is solved using the MUMPS direct solver \cite{amestoy2001mumps}.
	\item No initial-grid solver is employed during the iteration. We have observed no effect on convergence when the initial grid does not resolve the wave (as is the case throughout this work).
\end{itemize}
All experiments in this section were performed on \emph{Frontera}'s Cascade Lake (CLX) nodes at the Texas Advanced Computing Center \cite{stanzione2020frontera}. Timing statistics are neglected in this section, but will be provided in Section~\ref{sec:application}.

%=====================================================%
\subsection{Problem Setup}

Throughout this section we consider propagation of a Gaussian beam with waist-radius $0.1$ in a homogeneous unit cube domain $[0,1]^3$. Homogeneous impedance boundary conditions are imposed on all surfaces except near the origin where a Gaussian beam is injected through a manufactured impedance load. The solution is depicted in Fig.~\ref{fig:hp_beam}. 

%=====================================================%
\subsection{Direct Assembly vs.~Fine-Grid Restriction for Coarse-Grid Operators}
\label{sec:recomp}

As indicated in Section~\ref{sec:DPG-MG}, coarse-grid systems can be either directly assembled (or stored from previous meshes) or computed from fine-grid systems by applying the restriction operator. The two approaches are referred to as \emph{store} and \emph{restrict}, respectively. As will be demonstrated, the construction of coarse-grid systems has significant implications for the convergence of the DPG-MG solver.

\paragraph{Uniform $h$-refinements.}
%Insert error vs. error indicator here

\begin{figure}[!htb]
	\centering
	\begin{subfigure}{0.48\textwidth}
		\centering
		\includegraphics[width=\textwidth,trim={0pt 10pt 0pt 10pt},clip]
		{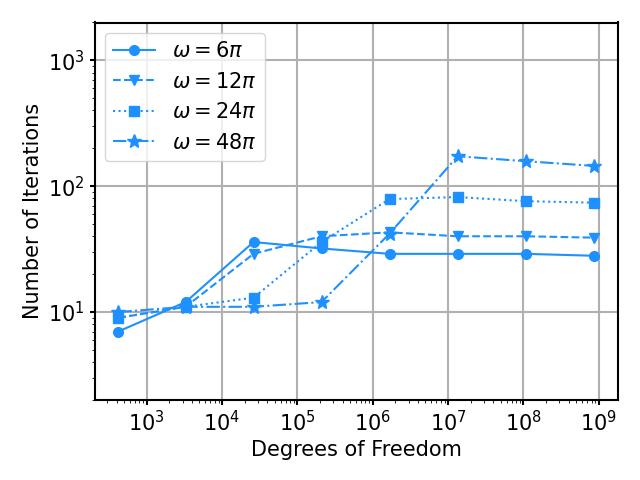}
		\caption{Coarse-grid system restricted from fine-grid system}
		\label{fig:h-a}
	\end{subfigure}
	\begin{subfigure}{0.48\textwidth}
		\centering
		\includegraphics[width=\textwidth,trim={0pt 10pt 0pt 10pt},clip]
		{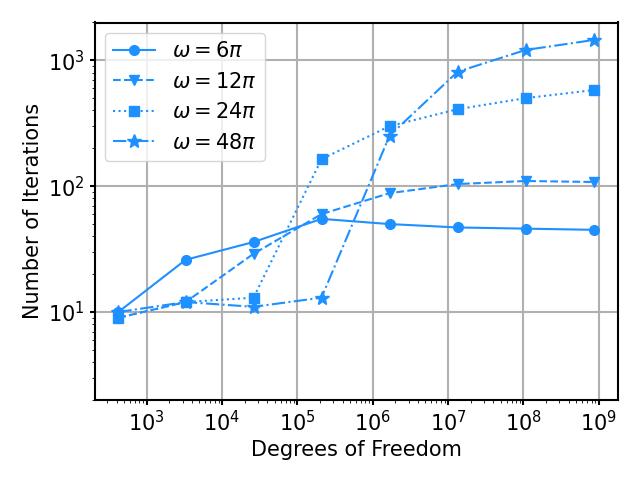}
		\caption{Coarse-grid system stored from previous meshes}
		\label{fig:h-b}
	\end{subfigure}
	\caption{Convergence of the DPG-MG solver with \emph{one smoothing step} applied to uniformly $h$-refined meshes. The solver convergence is $h$-robust and the iteration numbers are lower when using fine-grid restriction. The iterations until convergence depend linearly on the frequency $\omega$.}
	\label{fig:h}
\end{figure}

We begin by studying convergence of the DPG-MG solver under uniform $h$-refinements; i.e.~each subsequent grid is produced by a uniform $h$-refinement of the previous grid. The number of iterations required for convergence under each of the approaches, for a variety of frequencies, is reported in Fig.~\ref{fig:h}. Examining the results in Fig.~\ref{fig:h-a} (\emph{restrict}), it can be seen that the number of iterations increases roughly linearly with frequency but demonstrates clear $h$-robustness in the asymptotic regime. The increase in number of iterations with frequency is expected: meshes that cannot resolve the wave do not contribute to preconditioning the operator. Note that unlike multigrid preconditioners for the standard Galerkin method, which would diverge in this setting due to lack of discrete stability on the coarse grid, the DPG solver remains stable. Next, comparing Fig.~\ref{fig:h-a} (\emph{restrict}) and Fig.~\ref{fig:h-b} (\emph{store}), it can be seen that storing the coarse-grid system consistently resulted in a larger number of iterations than restricting; additionally, storing the coarse-grid system does not demonstrate $h$-robustness.

\paragraph{Uniform $h$-refinements; five smoothing steps.}

\begin{figure}[!htb]
	\centering
	\begin{subfigure}{0.48\textwidth}
		\centering
		\includegraphics[width=\textwidth,trim={0pt 10pt 0pt 10pt},clip]
		{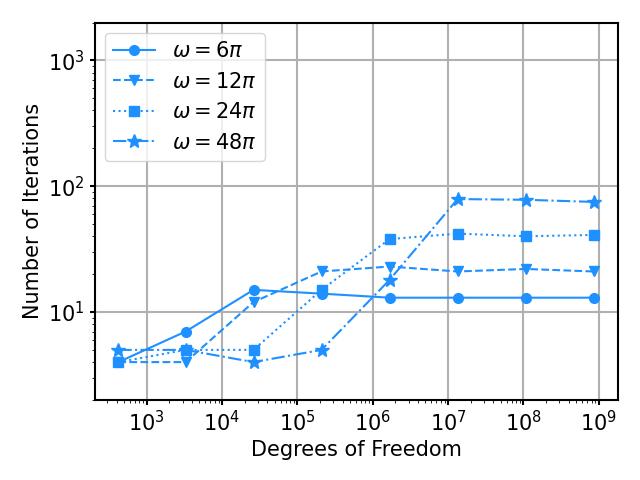}
		\caption{Coarse-grid system restricted from fine-grid system}
		\label{fig:h5-a}
	\end{subfigure}
	\begin{subfigure}{0.48\textwidth}
		\centering
		\includegraphics[width=\textwidth,trim={0pt 10pt 0pt 10pt},clip]
		{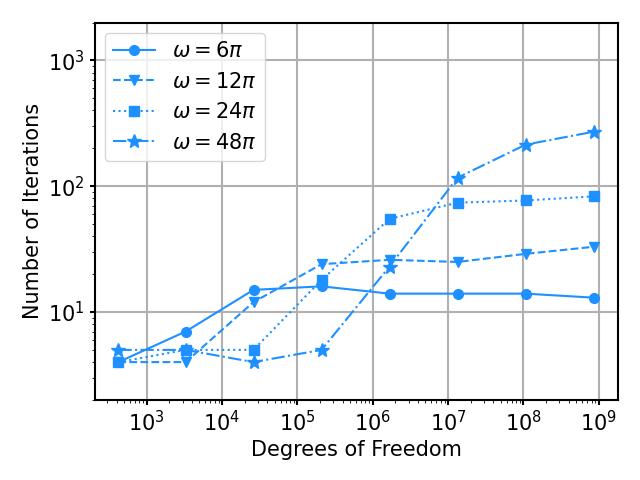}
		\caption{Coarse-grid system stored from previous meshes}
		\label{fig:h5-b}
	\end{subfigure}
	\caption{Convergence of the DPG-MG solver with \emph{five smoothing steps} applied to uniformly $h$-refined meshes. Doing additional smoothing tends to restore $h$-robustness in (b) to some extent; in (a), $h$-robustness is still observed, however using five smoothing steps per level only reduces the number of iterations until convergence by a factor of approximately two when compared to one smoothing step. Iterations depend linearly on the frequency $\omega$.}
	\label{fig:h5}
\end{figure}

In the original implementation of the DPG-MG solver, a relatively large number of smoothing iterations (typically between 5 and 10) were used in numerical experiments. For comparison, we repeat the previous study using five smoothing steps per iteration (V(5,5)-cycle); the results are depicted in Fig.~\ref{fig:h5}. Using a large number of smoothing steps tends to restore the $h$-robust convergence when the coarse-grid system is stored (Fig.~\ref{fig:h5-b}); however, the number of smoothing steps needed to attain $h$-robustness tends to increase with frequency. In particular, note that for the higher-frequency cases, the number of iterations until convergence is in fact lower when using one smoothing step with restriction from fine-grid systems (Fig.~\ref{fig:h-a}) than when using five smoothing steps with coarse-grid operators stored from previous meshes (Fig.~\ref{fig:h5-b}). Comparing Fig.~\ref{fig:h-a} and Fig.~\ref{fig:h5-a}, it can be seen that when the coarse-grid systems are defined via restriction, increasing the number of smoothing steps by a factor of five results in a decrease of the number of smoothing iterations by only a factor of two. We neglect an explicit study of convergence on the number of smoothing steps but qualitatively report that a single smoothing step per grid level is optimal for convergence in all of our numerical experiments to date.

\paragraph{Uniform $p$-refinements.}

\begin{figure}[!htb]
	\centering
	\begin{subfigure}{0.48\textwidth}
		\centering
		\includegraphics[width=\textwidth,trim={0pt 10pt 0pt 10pt},clip]
		{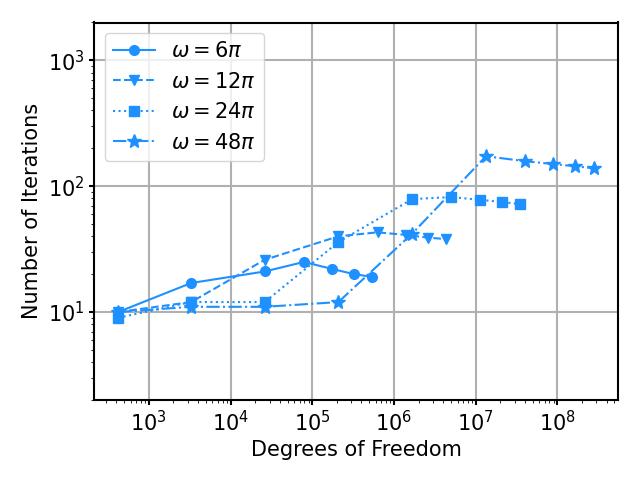}
		\caption{Coarse-grid system restricted from fine-grid system}
		\label{fig:p-a}
	\end{subfigure}
	\begin{subfigure}{0.48\textwidth}
		\centering
		\includegraphics[width=\textwidth,trim={0pt 10pt 0pt 10pt},clip]
		{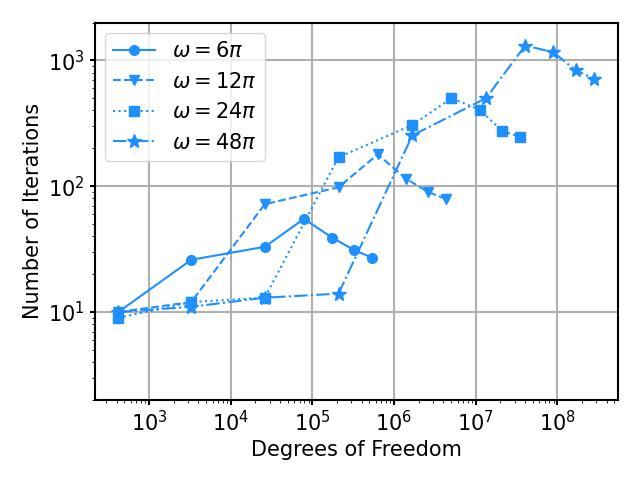}
		\caption{Coarse-grid system stored from previous meshes}
		\label{fig:p-b}
	\end{subfigure}
	\caption{Convergence of the DPG-MG solver applied to uniformly $hp$-refined meshes: grids are uniformly $h$-refined until two elements per wavelength, then uniformly $p$-refined. The solver convergence is $p$-robust and the iteration numbers are lower when using fine-grid restriction. The iterations until convergence depend linearly on the frequency $\omega$.}
	\label{fig:p}
\end{figure}

To investigate convergence of the DPG-MG solver under $p$-refinements, we first perform $h$-refinements until there are at least two elements per wavelength (satisfying the Nyquist criterion); then, the polynomial order of discretization $p$ is incremented on each subsequent grid. The results of this study are shown in Fig.~\ref{fig:p}, where it can be seen that both storing and restricting leads to $p$-robust convergence; however, restricting again requires far fewer iterations.

%=====================================================%
\subsection{\texorpdfstring{$hp$}{hp}-Adaptive Refinements}

\begin{figure}[!htb]
	\centering
	\includegraphics[width=\textwidth,trim={0pt 10pt 0pt 10pt},clip]
	{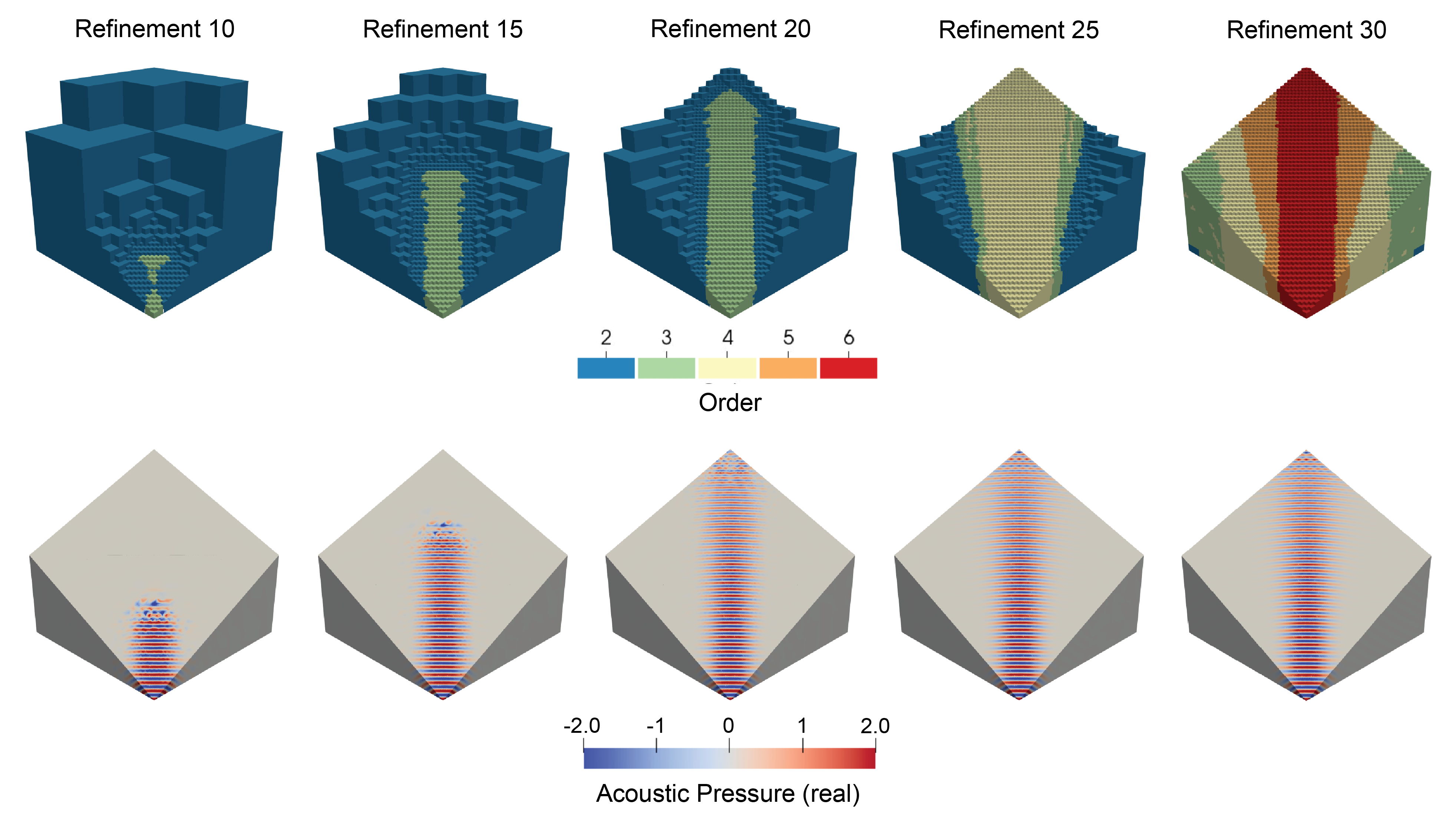}
	\caption{$hp$-adaptive propagation of a Gaussian beam in a cube domain; crinkle cut of $hp$-adaptive meshes (top) and surface cut of the real part of the acoustic pressure field (bottom). The $hp$-adaptive meshes have a ``sweeping'' structure as refinements first accumulate near the corner and then propagate into the domain; solutions on intermediate meshes are stable and partially resolve the wave.}
	\label{fig:hp_beam}
\end{figure}

We now consider $hp$-adaptive refinements, employing the D\"orfler marking strategy \cite{dorfler1996marking} to determine elements to be refined. Marked elements are $h$-refined until the maximum edge-length is less than one-half the wavelength, otherwise they are $p$-refined. We end refinements one mesh after no additional $h$-refinements are requested. As shown in Fig.~\ref{fig:hp_beam}, $hp$-adaptive refinements produce a series of meshes with a ``sweeping'' structure, i.e., they follow the direction of propagation of the beam.

\begin{figure}[!htb]
	\centering
	\begin{subfigure}{0.48\textwidth}
		\centering
		\includegraphics[width=\textwidth,trim={0pt 10pt 0pt 10pt},clip]
		{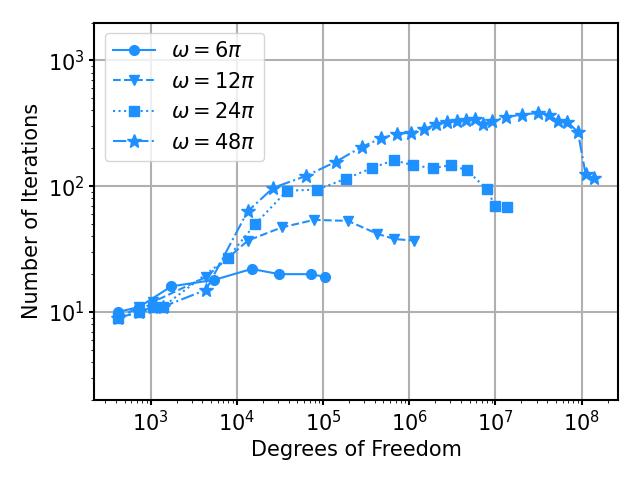}
		\caption{Number of iterations to $10^{-7}$ tolerance}
		\label{fig:hp-a}
	\end{subfigure}
	\begin{subfigure}{0.48\textwidth}
		\centering
		\includegraphics[width=\textwidth,trim={0pt 10pt 0pt 10pt},clip]
		{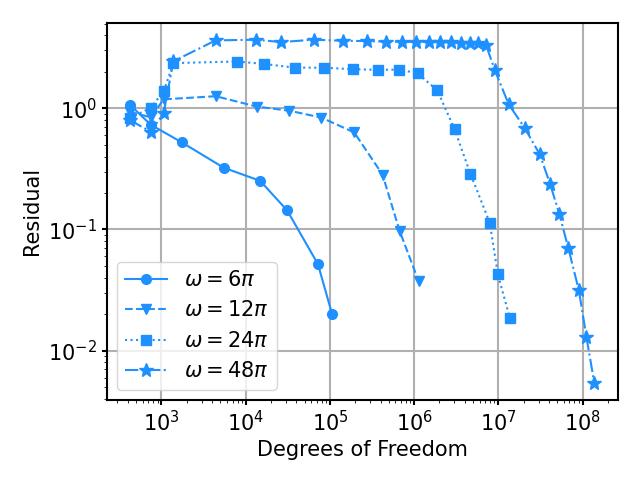}
		\caption{Convergence of DPG residual}
		\label{fig:hp-b}
	\end{subfigure}
	\caption{Convergence DPG-MG solver applied to $hp$-adaptive meshes. The number of iterations required for convergence again increases linearly with frequency ($\omega$) and the maximum number of iterations is higher than for uniform refinements (compare Fig.~\ref{fig:h-a}).}
	\label{fig:hp}
\end{figure}

In the case of uniform $h$-refinements, the observed linear increase in iterations with frequency is expected and is related to the inadequacy of coarse-space corrections when meshes are not sufficiently fine to resolve the wave. With $hp$-adaptive refinements, the behavior of the number of iterations until convergence with respect to frequency is less obvious since intermediate meshes are able to \emph{partially} resolve the wave. Indeed, we initially believed the ``sweeping'' structure of meshes helped to reduce the frequency dependence of convergence. The convergence study in Fig.~\ref{fig:hp} seems to indicate this is not the case; the number of iterations show a clear linear increase with frequency. However, note the maximum number of iterations required for convergence was consistently higher than for uniform refinements; this is unexpected since the adaptive case smoothes on each grid level, thus a much larger number of smoothing steps are performed overall. 

An initial hypothesis on the cause of the deteriorated convergence of the DPG-MG solver on $hp$-adaptive meshes implicates a similar phenomenon underlying the deteriorated convergence when coarse-grid systems are stored from previous meshes. In that case, roughly speaking, different scales between fine and coarse systems, imparted by element-wise Gram matrices, were ameliorated by constructing coarse-grid systems as restrictions of fine-grid systems, so that all systems inherit the fine-grid scale. Under $hp$-adaptive refinements, multiple element sizes---with various scales imparted by element-wise Gram matrices---are simultaneously present. We are undertaking a more rigorous investigation; however, note that in Fig.~\ref{fig:hp-a}, the number of iterations for convergence on the final meshes decreases considerably. Returning to Fig.~\ref{fig:hp_beam}, we can see that these final meshes are characterized exclusively by $p$-refinements with a consistent, uniform element size; from which coarse-grid systems are restricted.

%
%!TEX root = ../paper.tex
%

\section{Application -- Seismic Modeling}
\label{sec:application}

\begin{figure}[!htb]
	\centering
	\includegraphics[width=\textwidth,trim={0pt 0pt 30pt 50pt},clip]
	{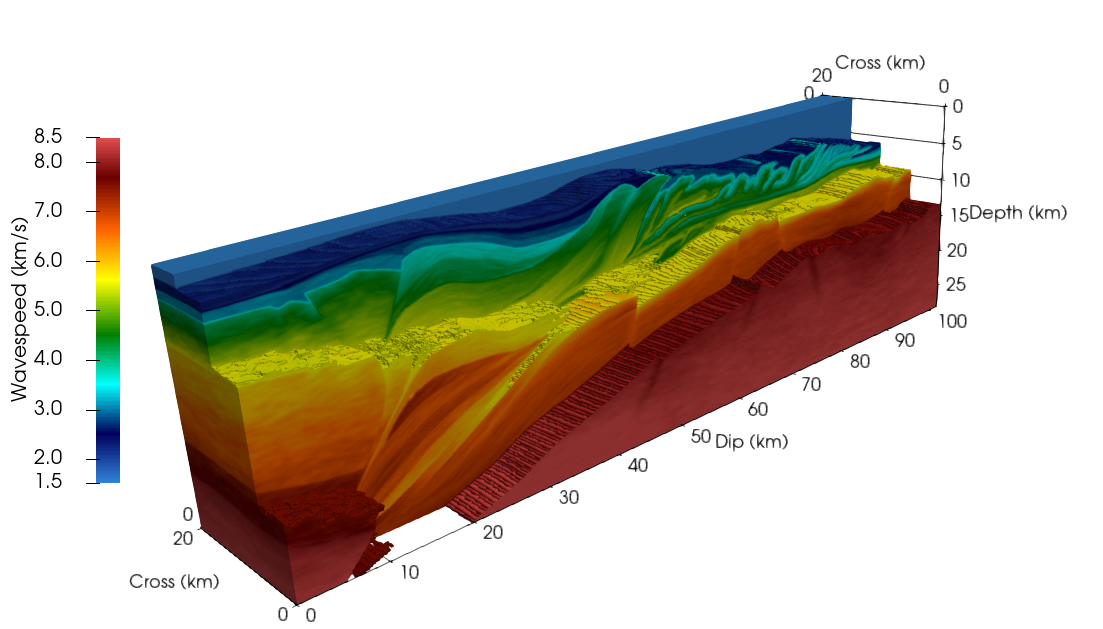}
	\caption{Cutaway of acoustic wavespeed for a section of the \gobench{} model, representing a subduction zone.}
	\label{fig:go-wavespeed}
\end{figure}

The observed linear increase of iterations with frequency $\omega$ implies a suboptimal $\mc{O}(N^{4/3})$ computational complexity of the DPG-MG solver in the preasymptotic regime; this is comparable to other methods including shifted Laplacian \cite{gander2015applying, sheikh2013convergence}, multilevel \cite{stolk2014multigrid, calandra2013improved, gopalakrishnan2004analysis}, and domain decomposition \cite{cai1992domain, kim2015optimized, stolk2013rapidly} methodologies, but is worse than the logarithmic increase observed for sweeping-type preconditioners \cite{engquist2011sweeping2,vion2014double} including source-transfer \cite{leng2020diagonal}, and others. Despite this linear increase, the DPG-MG solver is competitive for solving large-scale high-frequency wave propagation problems. To illustrate the performance and flexibility of the DPG-MG solver, we consider the \gobench{} benchmark \cite{go_3d_obs} from seismic modeling. The benchmark problem is set in a hexahedral domain with high-contrast heterogeneous structures representing a subduction zone, inspired by the geology of the Nankai Trough. Following \cite{tournier2022fefd}, we consider a $20 \times 102 \times 28.3$~km$^3$ section of the model. The wavespeed in this section is illustrated in Fig.~\ref{fig:go-wavespeed} and varies from $1\,500$~m/s to $8\,500$~m/s. Material data is specified on a uniform grid with spacing $100$~m, downsampled from the original $25$~m spacing of the \gobench{} model (which is rather large, occupying $132$~GB per parameter). The following computations were performed on $512$~\emph{Frontera} Cascade Lake (CLX) compute nodes ($28\,672$ cores) at the Texas Advanced Computing Center.

The problem is driven by a point source with a frequency of $3.75$~Hz located at $(10.0,12.5, 0.0)$, implemented as a tight Gaussian with standard deviation $\sigma=50$~m. For simplicity, we use an initial mesh consisting of $8 \times 42 \times 10$ hexahedral quadratic ($p=2$) elements with a total of $205\,757$ degrees of freedom. General unstructured meshes, fitted to high-contrast interfaces or adapted to wavespeed, could be used in conjunction with the DPG-MG solver with great effect; however, we illustrate that adaptive refinements, starting from an arbitrarily coarse initial mesh, can resolve complex problems without the need for hand-tuned or time-intensive meshing.

The D\"orfler marking strategy \cite{dorfler1996marking} is again used to mark elements for adaptive refinement. We perform seven initial $h$-adaptive refinements to resolve the region around the point source, followed by $hp$-adaptive refinements until the DPG residual is reduced by a factor of $3 \cdot 10^2$ (which coincided with the exhaustion of computer memory). The $hp$-adaptive strategy selects $h$-refinements when the maximum edge length of an element is less than one-half of the wavelength, except in regions of high contrast---which we define to be a greater than $10\%$ change in wavespeed over the element---where one-quarter of the wavelength is used. Otherwise, $p$-refinements are selected until a maximum order of $p=5$, after which $h$-refinements are again performed. Element order was limited to $p=5$ since higher-order elements can result in large smoothing patches that are expensive to store; we are pursuing a number of strategies to reduce patch storage, including a GPU implementation that recomputes smoothing patches during the solution. For applications that solve for many loads simultaneously, the cost of recomputing patches can be amortized over multiple loads, increasing the appeal of this approach. 

\begin{figure}[!htb]
	\centering
	\includegraphics[width=\textwidth]{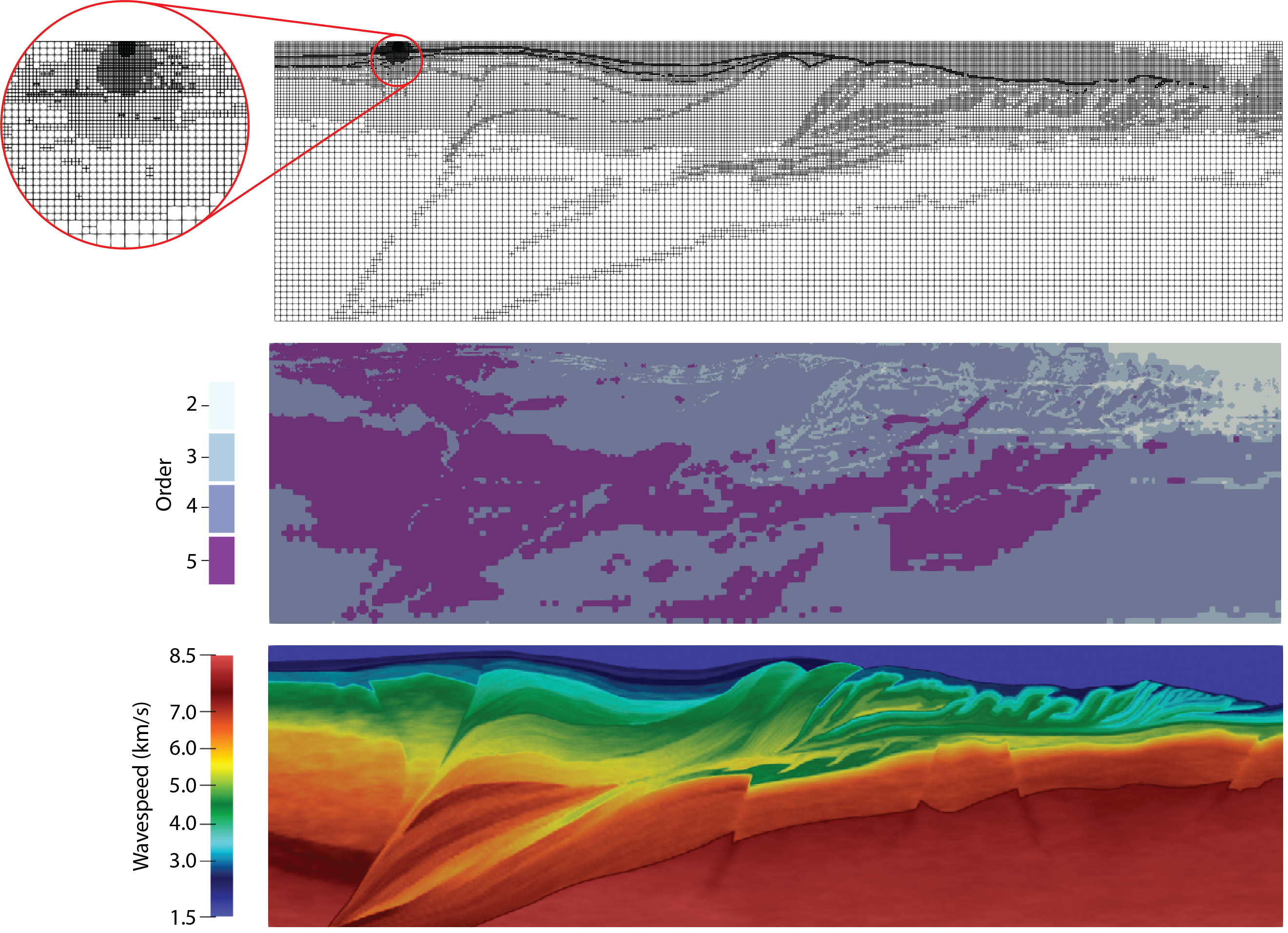}
	\caption{Cross-sections of the hexahedral mesh (top); element order (middle); and wavespeed (bottom), at a cross-wise distance of $10$~km. $h$-adaptivity is used to resolve high-contrast interfaces and the point source (enlarged to show detail) while $p$-adaptivity helps to mitigate the pollution error.}
	\label{fig:go-meshes}
\end{figure}

\begin{figure}[!htb]
	\centering
	\includegraphics[width=\textwidth]{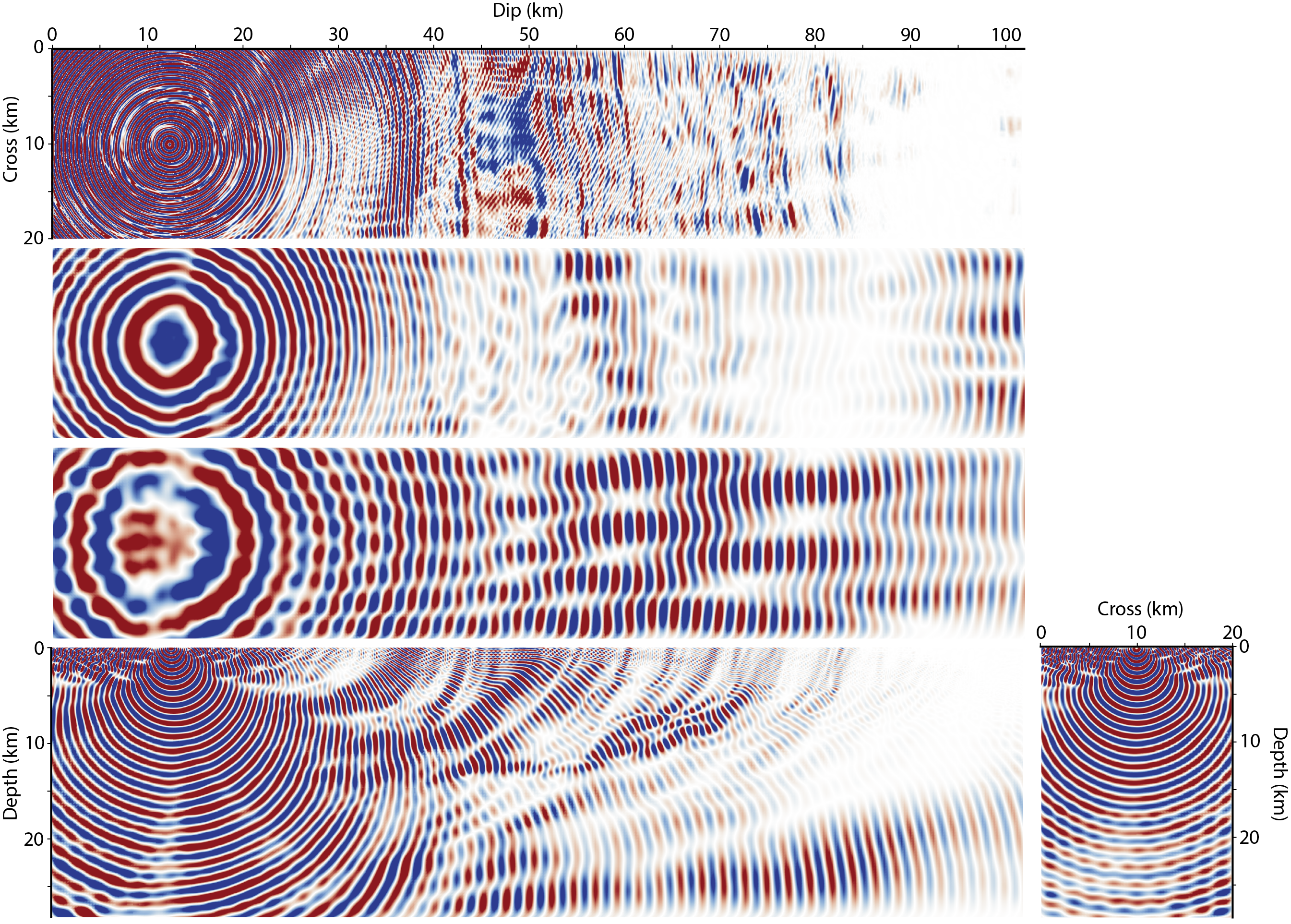}
	\caption{Acoustic pressure (real) for $3.75$~Hz frequency and a point source located at $(10,12.5,0)$. The following cross-sections are shown: from the top, $0.5$~km depth, $15$~km depth, $25$~km depth, and $10$~km cross-wise (left); and $12.5$~km dip-wise (right). The field is amplitude-compensated (scaled by the distance from the source) to enable visualization throughout the domain.}
	\label{fig:go-fields}
\end{figure}

In total, the mesh is refined 34 times, resulting in a final mesh (illustrated in Fig.~\ref{fig:go-meshes}) with over $6.3$ million elements and $1.9$ billion degrees of freedom. The upper right-hand corner of the mesh in Fig.~\ref{fig:go-meshes} is not fully refined and a few more adaptive steps would be needed to further refine the mesh in that region; however, looking ahead to Fig.~\ref{fig:go-fields}, it can be seen that the solution near this region has a relatively small amplitude. As alluded to previously, the refinement process was terminated due to a lack of memory; this is because the current implementation stores all $34$ adaptive meshes and associated coarse-grid operators (which are recomputed as restrictions of fine-grid operators after each refinement). Various amelioration strategies are possible, e.g.~refactoring the refinement tree to group refinements of similar depth and reduce the total number of grids.

% it feels strange to say the results are qualitatively similar when 
% a) the source is placed in a completely different place
% b) the solution of the other paper has a lot of errors we claim to resolve better

Cross-sections of the solution are illustrated in Fig.~\ref{fig:go-fields}. These results are qualitatively similar to those depicted in \cite{tournier2022fefd}, however we note the location of the point source between the two solutions differs slightly. A quantitative comparison of solution accuracy is deferred for a future work.

\begin{figure}[!htb]
	\centering
	\begin{subfigure}{0.48\textwidth}
		\centering
		\includegraphics[width=\textwidth,trim={0pt 10pt 0pt 10pt},clip]
		{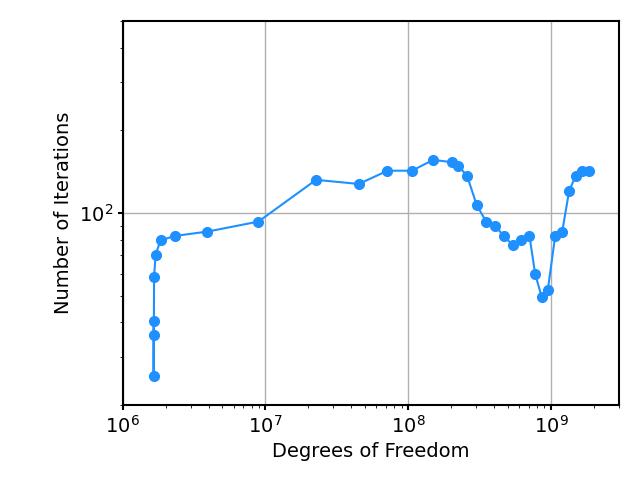}
		\caption{Number of iterations to convergence ($10^{-7}$ tolerance)}
		\label{fig:go-a}
	\end{subfigure}
	\begin{subfigure}{0.48\textwidth}
		\centering
		\includegraphics[width=\textwidth,trim={0pt 10pt 0pt 10pt},clip]
		{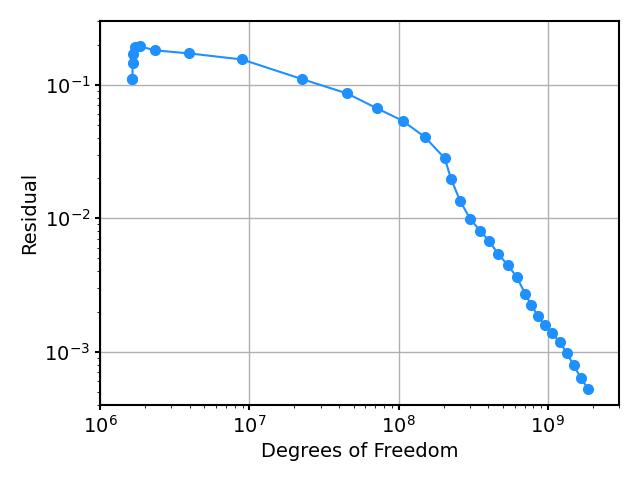}
		\caption{Convergence of DPG residual}
		\label{fig:go-b}
	\end{subfigure}
	\begin{subfigure}{0.48\textwidth}
		\centering
		\includegraphics[width=\textwidth,trim={0pt 10pt 0pt 10pt},clip]
		{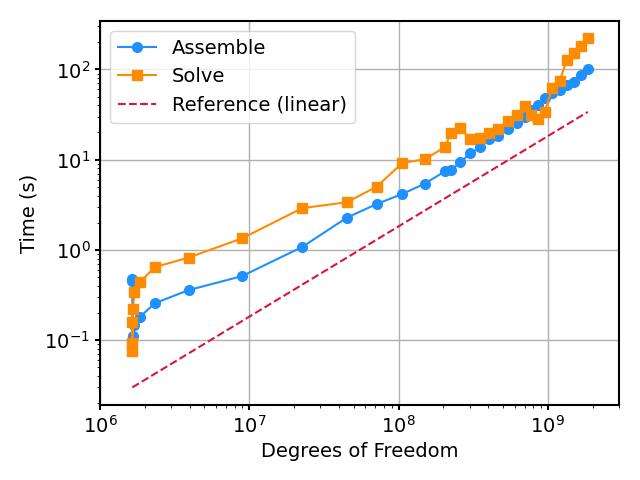}
		\caption{Assembly and solution time for the DPG system}
		\label{fig:go-c}
	\end{subfigure}
	\begin{subfigure}{0.48\textwidth}
		\centering
		\includegraphics[width=\textwidth,trim={0pt 10pt 0pt 10pt},clip]
		{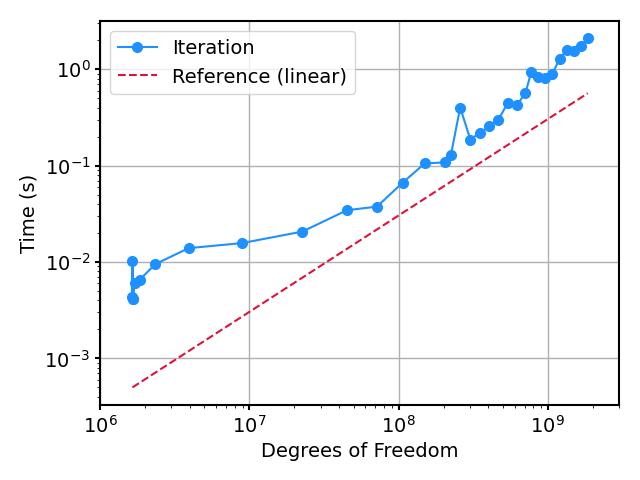}
		\caption{Total time \emph{per iteration}}
		\label{fig:go-d}
	\end{subfigure}	
	\caption{Convergence and timings for $hp$-adaptive solution of the \gobench{} model with the DPG-MG solver ($28\,672$ cores). The number of iterations is better controlled when solutions on previous grids are used to initialize subsequent grids. For a fixed frequency, the DPG-MG solver scales nearly linearly with respect to degrees of freedom; super-linear scaling is due to sub-geometric growth in DOFs.}
	\label{fig:go-stats}
\end{figure}

Convergence and timings for this example are shown in Fig.~\ref{fig:go-stats}. The number of iterations required for convergence (again using a tolerance of $10^{-7}$) on each grid is depicted in Fig.~\ref{fig:go-a}; however, note that in contrast to the convergence studies in Section~\ref{sec:studies}, here we use the solution on the previous meshes to initialize the solution on subsequent meshes. The effect of initializing with previous solutions becomes apparent in later iterations, when the solution is reasonably well resolved in much of the domain and further refinements result in fairly localized perturbations to the solution. Fig.~\ref{fig:go-b} shows that the DPG residual decreases early in the adaptive process, when compared to the $hp$-adaptive Gaussian beam problem in Fig.~\ref{fig:hp-b}; this is likely because the wave decays relatively quickly away from the point source.

Timings for the assembly and solution phase on each mesh are shown in Fig.~\ref{fig:go-stats}; the solution on the finest mesh was completed in $210$ seconds whereas the total runtime for the job (including assembly and solution on 34 meshes) was $3\,029$ seconds ($51$ minutes). Because the solution on different grids involve different numbers of iterations, the time \emph{per iteration} is depicted in Fig.~\ref{fig:go-d} to better gauge efficiency of the implementation. After an early preasymptotic regime, it can be seen that the time per iteration scales roughly linearly with respect to DOFs; however, some super-linearity is observed in the largest instances. Super-linearity of multigrid solvers is expected when the number of DOFs grows sub-geometrically between grid levels (as is the case late in the refinement process). Refactoring the refinement tree to reduce the number of grid levels could help mitigate super-linear scaling and significantly reduce the cost of $hp$-adaptive solver iterations. 

Finally, we note that the variations in time per iteration in Fig.~\ref{fig:go-d} are related to the significant challenge of load balancing on $hp$-adaptive meshes. Indeed, $hp$-adaptive refinements produce elements and smoothing patches with \emph{highly} disparate costs that must be accurately predicted then properly partitioned. We neglect definition of our load balancing strategy in this work, but we note that the multilevel approach employed here operates on a large number of relatively small elements and patches; this greatly simplifies the estimation of costs, can provide opportunities for more fine-grain parallelism, and is often more conducive to shared-memory parallelism than other methods including domain decomposition and sweeping-type preconditioners. Scalable adaptivity is a key differentiator of the DPG-MG solver which, to our knowledge, is novel among Helmholtz solvers.

\paragraph{Uniform refinements.}

We conclude this section by considering uniform $h$- and $p$-refinements for the \gobench{} benchmark. As remarked earlier, under the current implementation the DPG-MG solver stores and operates on all previous grids, which can become expensive when a small number of elements are refined between meshes. This added expense can be somewhat justified since adaptive meshes often attain a similar accuracy with a small fraction of the number of DOFs. Still, we are working toward reducing the additional expense of applying the solver to adaptively refined meshes. The following uniform refinement example is intended to provide a baseline for potential performance of adaptive refinements and to provide a more direct point of comparison to other Helmholtz solver implementations.

\begin{figure}[!htb]
	\centering
	\begin{subfigure}{0.48\textwidth}
		\centering
		\includegraphics[width=\textwidth,trim={0pt 10pt 0pt 10pt},clip]
		{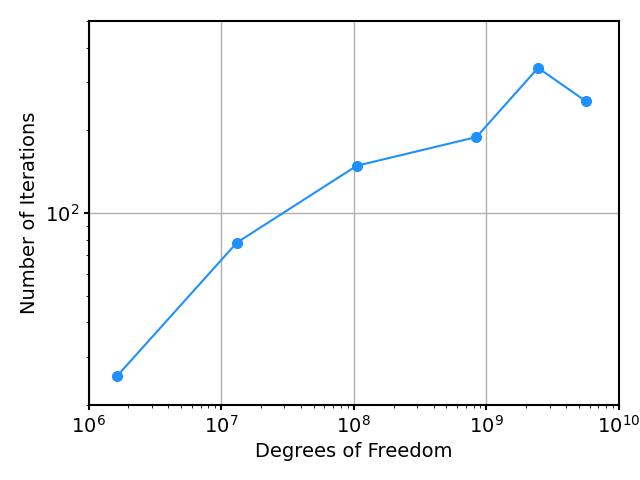}
		\caption{Number of iterations to convergence ($10^{-7}$ tolerance)}
		\label{fig:go-a-unif}
	\end{subfigure}
	\begin{subfigure}{0.48\textwidth}
		\centering
		\includegraphics[width=\textwidth,trim={0pt 10pt 0pt 10pt},clip]
		{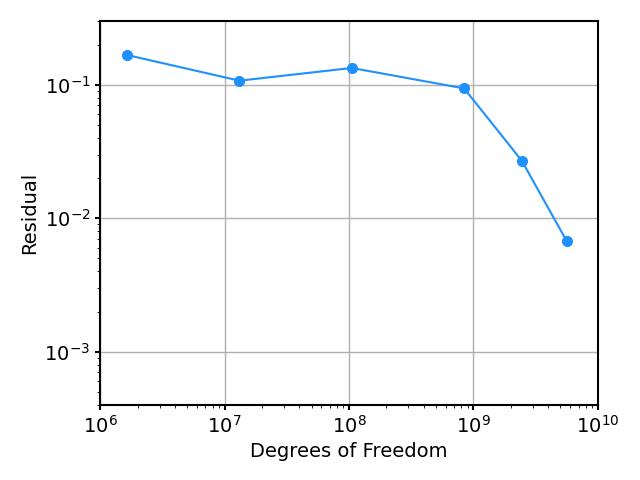}
		\caption{Convergence of DPG residual}
		\label{fig:go-b-unif}
	\end{subfigure}
	\begin{subfigure}{0.48\textwidth}
		\centering
		\includegraphics[width=\textwidth,trim={0pt 10pt 0pt 10pt},clip]
		{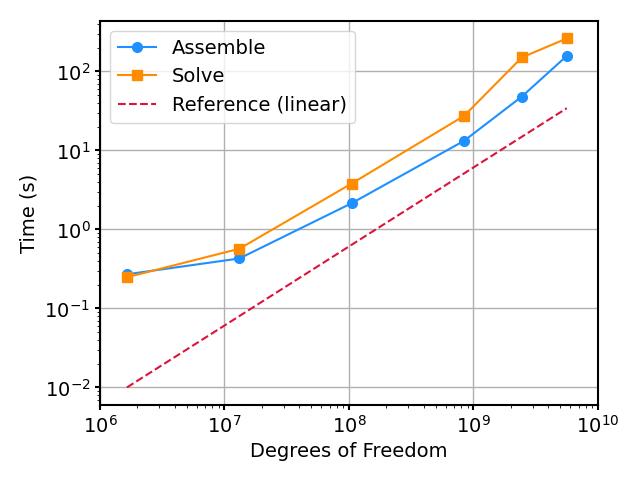}
		\caption{Assembly and solution time for the DPG system}
		\label{fig:go-c-unif}
	\end{subfigure}
	\begin{subfigure}{0.48\textwidth}
		\centering
		\includegraphics[width=\textwidth,trim={0pt 10pt 0pt 10pt},clip]
		{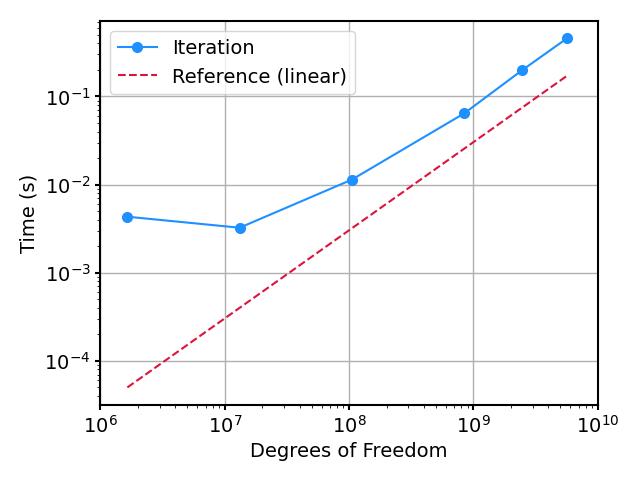}
		\caption{Total time \emph{per iteration}}
		\label{fig:go-d-unif}
	\end{subfigure}	
		\caption{Convergence and timings for uniform $h$- and $p$-refinements of the \gobench{} model with the DPG-MG solver ($28\,672$ cores). Nearly three times as many DOFs are required and the final DPG residual is ten times larger than the $hp$-adaptive case (compare Fig.~\ref{fig:go-b}); however, the time per iteration is three times smaller (compare Fig.~\ref{fig:go-d}).}
	\label{fig:go-stats-unif}
\end{figure}

We start from the same $8 \times 42 \times 10$ initial mesh of uniform order $p=2$ and perform four uniform $h$-refinements (corresponding to roughly two elements per wavelength in the water), followed by two uniform $p$-refinements. The final mesh has $16$ million elements and $5.6$ billion degrees of freedom. Convergence and timing information are provided in Fig.~\ref{fig:go-stats-unif}. Comparing Fig.~\ref{fig:go-b-unif} to Fig.~\ref{fig:go-b}, it can be seen that the uniform refinement setting requires nearly three times as many DOFs and only reaches a DPG residual ten times larger than for $hp$-adaptive refinements, illustrating the optimality of meshes produced with $hp$-adaptive refinements using the DPG error indicator. Still, comparing Fig.~\ref{fig:go-d-unif} and Fig.~\ref{fig:go-d}, it can also be seen that the time per iteration for uniformly refined meshes is three times smaller than for adaptive meshes, or nine times smaller when normalizing for the number of DOFs. A significant benefit may thus be attained by reducing the number of grid levels used during $hp$-adaptive refinements. Finally, we note that the number of iterations for uniform refinements (Fig.~\ref{fig:go-a-unif}) was higher than for the adaptive case (Fig.~\ref{fig:go-a-unif}); this is contrary to findings in Section~\ref{sec:studies} but seems to be the case here because adaptive refinements resolve the point source on early grids and thus typically begin iterations from a much lower residual than uniformly refined meshes.

%
%!TEX root = ../paper.tex
%

\section{Conclusion}
\label{sec:conclusion}

The DPG-MG solver leverages the unique properties of the DPG methodology including mesh-independent stability, a built-in error indicator, and Hermitian positive-definite discrete systems to enable robust, adaptable, and scalable solution of high-frequency wave propagation problems. When coarse-grid operators are constructed as restrictions of fine-grid operators, the DPG-MG solver demonstrates clear $h$- and $p$-robust convergence and a linear dependence with respect to wave frequency. A similar linear dependence on frequency was observed for $hp$-adaptive refinements. Despite the linear increase in number of iterations with respect to frequency, a scalable MPI/OpenMP implementation of the DPG-MG solver was demonstrated to be competitive for high-frequency wave propagation problems. In initial, moderate-scale tests, the DPG-MG solver was able to solve a challenging high-contrast seismic modeling benchmark (\gobench{}) with $1.9$ billion DOFs on $hp$-adaptive meshes and $5.6$ billion DOFs on uniformly refined meshes; a larger scale than any work we are currently aware of for high-frequency wave propagation in heterogeneous media. A significantly smaller DPG residual was achieved when $hp$-adaptive meshes were used, however we defer quantitative comparisons of accuracy to a later publication.

\paragraph{Future directions.}

Scalable implementation of the DPG-MG solver has motivated a number of promising research directions. First, as indicated in Section~\ref{sec:studies}, we intend to further investigate the deteriorated convergence rate of the DPG-MG solver when coarse-grid operators are stored and in the case of $hp$-adaptive meshes. Second, we are working to integrate the solver with automatic, fully-anisotropic $hp$-adaptivity \cite{hpbook,hpbook2}, where the DPG error indicator is first used to mark isotropic $hp$-refinements, an optimal set of anisotropic $h$- and $p$-refinements is then extracted from the isotropic refinements, and the remaining refinements are finally discarded. We anticipate the combination of the DPG-MG solver with automatic $hp$-adaptivity will be competitive for complex boundary layer problems and other problems with highly anisotropic features. We are continuously working to improve scaling and performance of the DPG-MG solver; near-term improvements include integration of GPUs for a memory-efficient implementation, refactorization of the refinement tree to reduce the number of grids under adaptive refinements, and implementation of a fully distributed data structure in $hp$3D (a fairly light-weight but replicated data structure currently limits scalability).

Finally, we are working to apply the scalable DPG-MG solver to challenging problems in science and engineering. For example, we have implemented an ultraweak time-harmonic Maxwell model of a tokamak device on an unstructured tetrahedral mesh with numerous reentrant corners and complex features; $hp$-adaptivity is expected to be particularly advantageous in this application. Additional applications of interest include extension of an optical fiber amplifier model \cite{henneking2021fiber,henneking2022parallel} to bent and complex cross-section fibers, and implementation of ultraweak elastic Helmholtz for seismic modeling and, ultimately, seismic inversion.

\FloatBarrier

%=====================================================%
% Bibliography
\printbibliography[heading=bibintoc]
\setcounter{section}{0}

\end{document}